\documentclass[11pt]{article}
\usepackage{amssymb}\usepackage{amsmath}
\newcommand{\beq}[1]{ 
 \begin{equation}\label{#1}}
\newcommand{\fr}[2]{{\textstyle \frac{#1}{#2} }}
\newcommand{\be}{\begin{eqnarray}}
\newcommand{\ee}{\end{eqnarray}}

\newcommand{\one}{{\mathfrak 1}}
\newcommand{\two}{{\mathfrak 2}}
\newcommand{\three}{{\mathfrak 3}}

\def\ot{\otimes}
\def\om{\omega}
\def\klein{\scriptscriptstyle}
\def\pis{\pi_s}

\def\PP{\mathbb P}
\def\EE{\mathbb E}

\def\SJ{\mathsf J}
\def\SG{\mathsf G}
\def\SH{\mathsf H}
\def\SX{\mathsf X}
\def\SO{\mathsf O}
\def\SR{\mathsf R}
\def\SP{\mathsf P}
\newcommand{\rf}[1]{(\ref{#1})}
\newcommand{\+}{\,{+}\,}
\renewcommand{\=}{\,{=}\,}
\newcommand{\rfeq}[1]{\stackrel{(\ref{#1})}{=}}
\newcommand{\FRW}[6]{\Bigl\{\begin{smallmatrix}%
 \textstyle #1\, & \textstyle #3\, & \textstyle #5 \\[2pt]%
 \textstyle #2\, & \textstyle #4\, & \textstyle #6 %
 \end{smallmatrix}\Bigr\}}

\newtheorem{prop}{Proposition}
\newtheorem{lem}{Lemma}
\newcounter{num}
\begin{document}
\strut\hfill March  2010\\
\vspace*{3mm}
\begin{center}
{\large\bf
 On constant $U_q(sl_2)$--invariant R--matrices }
\\ [3mm]
{\sc A. G. Bytsko} \\ [2mm]
{ Steklov Mathematics Institute \\
 Fontanka 27, 191023, St.-Petersburg, Russia }

\end{center}
\vspace{1mm}
\begin{abstract}
\noindent
 The spectral resolution of a $U_q(sl_2)$--invariant
 solution $R$ of the constant Yang--Baxter equation
 in the braid group form is considered. It is shown
 that, if the two highest coefficients in this
 resolution are not equal, then $R$ is either the Drinfeld
 R--matrix or its inverse.

\end{abstract}
\subsection*{\S 1. Introduction}\label{AUR}   

Recall that the algebra $U_q(sl_2)$ is generated by the generators
$X_+$, $X_-$, $q^{H}$, $q^{-H}$ satisfying the relations
\cite{KR2}
\beq{Uq}
 [X^+,X^-] = \fr{q^{2H}-q^{-2H}}{q-q^{-1}} \,, \qquad
 q^{H} X^\pm = q^{\pm 1} \, X^\pm \, q^{H} \,,\qquad
 q^{\pm H} q^{\mp H} = 1\,.
\end{equation}
The homomorphism $\Delta$ which is defined on the generators
as follows
\beq{DelEFK}
 \Delta(X^\pm) = X^\pm \ot q^{-H} + q^{H}\ot X^\pm \,,
 \qquad \Delta(q^{\pm H}) = q^{\pm H} \ot q^{\pm H} \,,
\end{equation}
turns $U_q(sl_2)$ into a bialgebra (moreover, a Hopf
algebra \cite{Sk}).

We will consider the standard finite dimensional representation
$\pis$ of the algebra $U_q(sl_2)$ in which the generators act on
the basis vectors $\om_{k}$ of a module $V_s$
($\mathrm{dim}\, V_s \= (2s{+}1)$,
$2s \,{\in}\, {\mathbb N}$) as follows
\beq{EFK}
 \pis(X^\pm)\,\om_k = \sqrt{[s \mp k][s \pm k + 1]} \ \om_{k \pm 1} \,,
\qquad
 \pis(q^{\pm H})\,\om_k = q^{\pm k} \, \om_{k} \,,
\end{equation}
where $[t] \equiv (q^t \,{-}\, q^{-t})/(q \,{-}\, q^{-1})$
and $k \= {-}s,{-}s {+} 1,{\ldots},s$.

The universal R--matrices for the algebra \rf{Uq}--\rf{DelEFK}
are given by \cite{D1}
\beq{Ru}
  R^\pm = q^{\pm H \ot H} \, \sum_{n =0}^\infty
 \frac{ q^{\pm\frac{1}{2}(n^2-n)} }{\prod_{k=1}^n [k]_q}
 \bigl( \pm (q\,{-}\,q^{-1}) X^\mp \ot X^\pm \bigr)^n  \,
 q^{\pm H \ot H} \,.
\end{equation}

Let $\PP$ denote the operator which permutes the tensor
components in $U_q(sl_2)^{\otimes 2}$. Then the operator
$R \,{\equiv}\, \PP \, R^+ \= (R^-)^{-1} \, \PP$
satisfies the Yang--Baxter equation in the braid group form:
\beq{YB}
   R_{\one\two} \, R_{\two\three} \, R_{\one\two} =
   R_{\two\three} \, R_{\one\two} \, R_{\two\three} \,.
\end{equation}
The spectral resolution of $R$ in the representation $\pi_s$
is given by \cite{KR1}
\beq{spec}
 \SR \equiv  \pis^{\ot 2} \, (R) =
  \sum_{k=0}^{2s} \ \xi_k \, \SP^{2s-k} \,,
\end{equation}
where $\SP^j$ stands for the projector onto the irreducible
submodule $V_j$ in~$V_s^{\ot 2} = {\oplus}_{j=0}^{2s} V_j$.
Here and below we use the following notations
\beq{rho}
 \xi_k \equiv (-1)^{k} \,  q^{\rho(2s-k) -2\rho(s)} \,, \qquad
 \rho(t)\equiv t(t \+1) \,.
\end{equation}

Consider an $U_q(sl_2)$--invariant solution $R'$ of the
Yang--Baxter equation~\rf{YB}. Its spectral resolution in the
representation $\pi_s$ is given by
\beq{R'}
\SR' \,{\equiv}\, \pi_s^{\ot 2} (R') =
\sum_{k=0}^{2s} r_k \SP^{2s-k} \,,
\end{equation}
where $r_0 \,{\neq} \,0$ by Lemma~6 в \cite{B2}, which
applies to the case $q \,{\neq}\,1$ as well.
We will prove the following statement.

\begin{prop}\label{Hc0}
If the spectral resolution \rf{R'} has $r_1 \neq r_0$,
then $\SR'$ coincides with either $\SR$
or with $\SR^{- 1}$ up to normalization.
\end{prop}

This statement is a $q$--analogue of the second part of
Proposition~1 in \cite{B2}, where $sl_2$--invariant
solutions of the Yang--Baxter equation were considered.
Note that the limit $q \to 1$ is degenerate in the sense that
both operators $\SR$ and $\SR^{-1}$ turn into the
permutation operator~$\PP$.

\subsection*{\S 2. Reduction on the subspace $W^{(s)}_n$}   

Let us recall the method of analyzing $U_q(sl_2)$--invariant
solutions of the Yang--Baxter equation developed in~\cite{B1}.
Let $\lfloor t \rfloor$ denote the entire part of~$t$.
The subspace $W^{(s)}_n \subset V_s^{\otimes 3}$ for
\hbox{$n=0,1,\ldots,\lfloor 3s \rfloor$} is defined as a span
of the highest weight vectors of weight $(3s \,{-}\, n)$, that is
\beq{Wn}
  W^{(s)}_n = \{ \, \psi \in V_s^{\otimes 3} \quad \bigm| \quad
  \SX^+_{\one\two\three} \psi =0 \,, \quad
  q^{\SH_{\one\two\three}} \psi = q^{3s - n} \psi \, \} \,.
\end{equation}
Here and below for $O \in U_q(sl_2)$ we use the notation:
$\SO_{\one\two\three} = \pis^{\ot 3}\bigl((\Delta \ot id)\Delta(O)\bigr)$.

Since
$[\SX^\pm_{\one\two\three},\SR_{\one\two}] \=
[\SX^\pm_{\one\two\three},\SR_{\two\three}] \=0$,
then $W^{(s)}_n$ is an invariant subspace for $\SR_{\one\two}$ and
$\SR_{\two\three}$, and so we can consider reductions
of these operators onto $W^{(s)}_n$.
We can choose a basis of $W^{(s)}_n$ in which the operator
\hbox{$\SR_{\one\two}\!\!\bigm|_{W^{(s)}_n}$}
is represented by a diagonal matrix $D^{(n)}_0$ of the following form:
\beq{Dkk}
 \bigl(D^{(n)}_0\bigr)_{k k'} = \delta_{k k'} \, \xi_k \,,
\end{equation}
here $0 \leq k \leq n$ for $0 \leq n \leq 2s$ and
$(n \,{-}\, 2s) \leq k \leq (4s \,{-}\, n)$ for
$2s \leq n \leq \lfloor 3s \rfloor$.

In the same basis, the operator
$\SR_{\two\three}\!\!\bigm|_{W^{(s)}_n}$ is represented by
the following matrix
\beq{Dhat}
 \hat{D}^{(n)}_0 = A^{(s,n)} \ D^{(n)}_0 \ A^{(s,n)} \,,
\end{equation}
where $A^{(s,n)}$ is a matrix with the following
properties \cite{B1}: it is symmetric, orthogonal,
equal to its inverse, and self--dual in  $q$:
\beq{A3}
 A^{(s,n)} =  \bigl( A^{(s,n)} \bigr)^t =
  \bigl( A^{(s,n)} \bigr)^{-1} \,, \qquad
 A^{(s,n)}_q = A^{(s,n)}_{q^{-1}}\,.
\end{equation}
Its entries are expressed in terms of the 6--$j$ symbols
of the algebra $U_q(sl_2)$ as follows:
\beq{A1}
 A^{(s,n)}_{k k^{\prime} }  =
 (-1)^{2s-n} \sqrt{[4s- 2k + 1]_q [4s - 2k' + 1]_q} \
 \FRW{s}{s}{s}{3s - n}{2s - k}{2s - k'}_q \,.
\end{equation}

The statement that the Yang--Baxter equation \rf{YB} holds
when it is reduced onto the subspace ~$W^{(s)}_n$
is equivalent to the following equality
\beq{RYB0}
 \bigl( D^{(n)}_0  \, A^{(s,n)} \bigr)^3 =
 \bigl( A^{(s,n)} \, D^{(n)}_0 \bigr)^3 \,.
\end{equation}
Actually, however, a stronger statement holds:
the r.h.s. and the l.h.s. of \rf{RYB0} are equal up
to a multiplicative constant to the identity
operator on $W^{(s)}_n$.
This follows from the following statement (which
is a q--analogue of Lemma ~3 in \cite{B2}):

\begin{lem}\label{PM}
For all \hbox{$n=0,\ldots,{\lfloor 3s \rfloor}$},
the following relation holds:
\beq{AD0}
  A^{(s,n)} \,  D^{(n)}_0  \, A^{(s,n)} = \theta_n \,
 \bigl( D^{(n)}_0 \bigr)^{-1}\, A^{(s,n)} \,
 \bigl( D^{(n)}_0 \bigr)^{-1} \,,
\end{equation}
where $\theta_n \equiv (-1)^n \, q^{\rho(3s-n) -3\rho(s)}$.
\end{lem}

The proof of this and other lemmas is given in the Appendix.

The statement of Lemma~\ref{PM} can be written in the
following form:
\beq{ada}
 \bigl( \SR_{\one\two} \, \SR_{\two\three} \,
    \SR_{\one\two} \bigr) \!\!\bigm|_{W_n^{(s)}} =
 \bigl( \SR_{\two\three} \, \SR_{\one\two} \,
    \SR_{\two\three} \bigr) \!\!\bigm|_{W_n^{(s)}}
 = \theta_n \, A^{(s,n)} \,.
\end{equation}
For $q=1$ this relation turns into
$\bigl( \PP_{\one\three} \bigr) \!\!\bigm|_{W_n^{(s)}} =
 (-1)^n A^{(s,n)}$.

{}From \rf{ada} and \rf{A3} it follows that
\beq{ada2}
 \bigl( (\SR_{\one\two} \, \SR_{\two\three} )^3
      \bigr) \!\bigm|_{W_n^{(s)}} =
 \bigl( (\SR_{\two\three} \, \SR_{\one\two} )^3
      \bigr) \!\bigm|_{W_n^{(s)}}
 =   q^{2\rho(3s-n) -6\rho(s)}  \,.
\end{equation}

Let us note that
\begin{eqnarray}
{}&&  \bigl( \SR_{\one\two} \, \SR_{\two\three} \,
    \SR_{\one\two} \bigr)^2
 = \bigl( \SR_{\two\three} \, \SR_{\one\two} \,
    \SR_{\two\three} \bigr)^2
 = \bigl( \SR_{\one\two} \, \SR_{\two\three} \bigr)^3
 = \bigl( \SR_{\two\three} \, \SR_{\one\two} \bigr)^3 \\
\label{ada3}
{}&& = \pi_s^{\ot 3} \Bigl( \bigl(
    R^-_{\one\two} R^-_{\one\three} R^-_{\two\three} \bigr)^{-1}
 \bigl( R^+_{\one\two} R^+_{\one\three} R^+_{\two\three} \bigr) \Bigr)
 = \pi_s^{\ot 3} \bigl( \chi_\one \chi_\two \chi_\three \,
    \Delta^{(2)} (\chi^{-1})  \bigr) \,,
\end{eqnarray}
where the element $\chi$ is constructed in the following way:
write the $R$--matrix \rf{Ru} as
$R^+ \= \sum_a r^{\klein(1)}_a \,{\ot}\, r^{\klein(2)}_a$,
and let ${\cal S}$ stand for the antipode operation, then
$\chi \= q^{2H}  \bigl(\sum_a {\cal S}(r^{\klein(2)}_a)
 r^{\klein(1)}_a \bigr)$.
It is known \cite{D2} that the element $\chi$ is central,
$\pi_s(\chi)=q^{-2\rho(s)}$, and
$\chi_1 \chi_2 \, \Delta(\chi^{-1}) =
 \bigl(R^- \bigr)^{-1}  R^+$.
The last relation allows us to derive the last equality in
\rf{ada3} (and its generalization for $\Delta^{(N)}(\chi^{-1})$,
see the proof of Lemma~1 in~\cite{B3}).
Thus, relation \rf{ada} can be regarded as the definition of
a certain square root of the operator given by the r.h.s. of~\rf{ada3}.

\subsection*{\S 3. Yang--Baxter equation on
  $W_n^{(s)}$ }   

We will prove Proposition~1 using the following
statement (a $q$--analogue of Lemma~4 in~\cite{B2}.

\def\ovm{\overline{m}}
\begin{lem}\label{PPrel1}
Let $0 \,{\leq}\, \ovm \,{\leq}\, n \,{\leq}\, 2s$,
where $\ovm \,{\equiv}\, (2s \,{-}\,m)$.
The reductions of the operators
$\SP^{m}_{\one\two}$, $\SP^{m}_{\two\three}$,
and $\SR_{\one\two}^{\pm 1}$, $\SR_{\two\three}^{\pm 1}$ on~$W_n^{(s)}$
satisfy the following relations
\be
\label{pp3}
  \SR_{l} \, \SR_{l'} \, \SR_{l} = \SR_{l'} \, \SR_{l} \, \SR_{l'} , &&
 \SP^{m}_{l} \, \SP^{m}_{l'} \, \SP^{m}_{l}
    = \eta_{n,\ovm}^2 \, \SP^{m}_{l} , \\
\label{pp4}
 \SP^{m}_{l} \, \SR_{l'}^{\pm 1} \, \SP^{m}_{l} =
    (\theta_n \xi_{\ovm}^{-2})^{\pm 1}
    \eta_{n,\ovm} \, \SP^{m}_{l} , &&
  \SR_{l}^{\pm 1} \, \SP^{m}_{l'} \, \SR_{l}^{\pm 1} =
    (\theta_n \xi_{\ovm}^{-1})^{\pm 2} \,
    \SR_{l'}^{\mp 1} \, \SP^{m}_{l} \, \SR_{l'}^{\mp 1} ,  \\ [1mm]
\label{pp5}
 \SP^{m}_{l} \, \SP^{m}_{l'} \, \SR_{l}^{\pm 1} =
    (\theta_n \xi_{\ovm}^{-1})^{\pm 1} \eta_{n,\ovm} \,
    \SP^{m}_{l} \, \SR_{l'}^{\mp 1} , &&
 \SR_{l}^{\pm 1} \, \SP^{m}_{l'} \, \SP^{m}_{l} =
 (\theta_n \xi_{\ovm}^{-1})^{\pm 1} \eta_{n,\ovm} \,
    \SR_{l'}^{\mp 1} \SP^{m}_{l} ,
\ee
where $l={\scriptstyle \{12\}}$, $l'={\scriptstyle \{23\}}$ or
$l={\scriptstyle \{23\}}$, $l'={\scriptstyle \{12\}}$,  and\,
$\eta_{n,\ovm} = A^{{(s,n)}}_{\ovm,\ovm}$.
\end{lem}

Let us remark that not all relations in Lemma~\ref{PPrel1} are
independent. For instance, the second relation in \rf{pp4}
follows from \rf{pp5}; the first relation in \rf{pp4} and
the second relation in \rf{pp3} can be derived from each
other with the help of~\rf{pp5}.

Let us remark also that, for $q \= 1$, the operators $\SR^{\pm 1}$
coincide with the permutation operator $\PP$, and relations
\rf{pp3}--\rf{pp5} become the relations of the Brauer algebra
\cite{Br} (taking into account the additional relation
$\PP^2 \= \EE$, where $\EE$ is the identity operator).
For $q \,{\neq}\, 1$, the reductions of the operators $\SR^{\pm 1}$
onto $W_1^{(s)}$ can be represented as linear combinations of
$\SP^m$ and the identity operator $\EE$. As a consequence,
relations \rf{pp3}--\rf{pp5} for $n \= 1$ can be derived from
the second relation in \rf{pp3}, which is the defining relation
for the Temperley--Lieb algebra~\cite{TL}.
For $n \,{\geq}\, 2$,  relations \rf{pp3}--\rf{pp5} are the
relations that hold in the Birman--Wenzl--Murakami algebra~\cite{BW,Mu}.
However, in this algebra an additional relation must also hold,
which in our case holds only for $n \= 2$ (the operator
$\SR^{-1}$ being reduced onto $W_2^{(s)}$ can be represented as
a linear combination of the operators $\SR$, $\SP^m$, and $\EE$).

Returning to consideration of the spectral resolution \rf{R'},
let us note that without a loss of generality we can set
$r_0 \,{=} \,\xi_0$.
Then $\SR'$ can be represented in the following form:
\beq{Rans}
 \SR' = \SR + g\, \SP^{2s-n} + \ldots \,,
\end{equation}
where $n\,{\geq} \,1$ and $\ldots$ stands for the
sum involving projectors of ranks smaller than the
rank of $\SP^{2s-n}$.

Substitute the ansatz \rf{Rans} in the Yang--Baxter
equation and consider its reduction onto $W_n^{(s)}$
for $n \leq 2s$.
With the help of relations of Lemma~\ref{PPrel1}, it can be
verified that the Yang--Baxter equation for
$\SR'\!\bigm|_{W_n^{(s)}}$  is equivalent to the
following matrix equation
\begin{equation}\label{FGH}
 g \, \SJ +
 (\theta_n  \xi_{n}^{-2}  \eta_{n,n} \, g^2
    + \eta_{n,n}^2 \, g^3) \, \SG  +
  (\theta_n  \xi_{n}^{-1}  \eta_{n,n} \, g^2) \, \SH = 0 \,,
\end{equation}
where
\begin{equation*}
\begin{aligned}
\SG &= \bigl( \SP^{2s-n}_{\one\two} -
    \SP^{2s-n}_{\two\three} \bigr) \!\bigm|_{W_n^{(s)}}  =
    \pi^{(n)} - A^{(s,n)} \pi^{(n)} A^{(s,n)} ,\\
\SJ &= \bigl( \SR_{\one\two} \, \SP^{2s-n}_{\two\three} \, \SR_{\one\two}
 - \SR_{\two\three} \, \SP^{2s-n}_{\one\two} \, \SR_{\two\three}
 \bigr) \!\bigm|_{W_n^{(s)}}  \\
{} & = D_0^{(n)} A^{(s,n)} \pi^{(n)} A^{(s,n)} D_0^{(n)} -
  \theta_n^2  \xi_{n}^{-2}
 (D_0^{(n)})^{-1} A^{(s,n)} \pi^{(n)} A^{(s,n)} (D_0^{(n)})^{-1} , \\
\SH &= \bigl( \SP^{2s-n}_{\one\two} \, \SR^{-1}_{\two\three} +
 \SR^{-1}_{\two\three} \, \SP^{2s-n}_{\one\two} -
 \SP^{2s-n}_{\two\three} \, \SR^{-1}_{\one\two} -
 \SR^{-1}_{\one\two} \, \SP^{2s-n}_{\two\three}
    \bigr) \!\bigm|_{W_n^{(s)}}   \\
{}& =  \theta_n^{-1}  \xi_{n}
 \, \bigl(\pi^{(n)} A^{(s,n)} D_0^{(n)} +
    D_0^{(n)} A^{(s,n)} \pi^{(n)} \bigr) \\
{} &\ \ \ - A^{(s,n)} \pi^{(n)} A^{(s,n)} (D_0^{(n)})^{-1}  -
    (D_0^{(n)})^{-1} A^{(s,n)} \pi^{(n)} A^{(s,n)} .
\end{aligned}
\end{equation*}
Here $\pi^{(n)}$  is a matrix such that
$(\pi^{(n)})_{kk'} = \delta_{kn}\delta_{k'n}$.

\begin{lem}\label{LD}
i) For $n = 1$, the following relations hold:
\beq{JGH1}
 \SJ = (\theta_1^2 \xi_0^{-2} \xi_1^{-2} - \xi_0^2) \,\SG =
 (q^{4s(s-1)} \,{-}\, q^{4s^2})\, \SG \,,\qquad
 \SH = 2 \xi_0^{-1} \, \SG = 2 q^{-2s^2} \SG \,.
\end{equation}
ii) For $n = 2$,  the matrices $\SJ$ and $\SG$ are linearly
independent, and the following relation holds:
\beq{JGH2}
 \xi_0 \xi_1\, \SH = (\xi_0 + \xi_1) \, \SG +
    (\xi_0 + \xi_1)^{-1} \, \SJ\,.
\end{equation}
iii) For $n \geq 3$, the matrices $\SJ$, $\SG$, $\SH$
are linearly independent, and $\SJ \,{\neq}\, 0$.
\end{lem}

Substituting relations \rf{JGH1} in \rf{FGH}, we infer that,
for $n \= 1$, the coefficient $g$ must be a root of the
following equation:
\begin{equation*}
 \eta_{1,1}^2 \, g^3 +
    \eta_{1,1} \theta_1 \xi_1^{-1} (\xi_1^{-1}+2\xi_0^{-1}) \, g^2
    +(\theta_1^2 \xi_0^{-2} \xi_1^{-2} - \xi_0^2) \, g =0 .
\end{equation*}
Hence, taking into account that
$\eta_{1,1} \={-}\,(q^{2s} \+ q^{-2s})^{-1}$,
we find that, for $n \= 1$, the coefficient $g$ can take one of the
following values:
$g \=0$, $g \= q^{2s(s-2)}(1 \,{-}\, q^{8s})$,
$g \= q^{2s(s-2)}(1 \+ q^{4s})$.
In the first and second cases, the spectral resolution of $\SR'$
coincides in the two highest orders with that of $\SR$ and
$q^{4s^2} \SR^{-1}$, respectively.
In the third case, we have $r_1 \= r_0$.

For $n \= 2$, substitute relations \rf{JGH2} in \rf{FGH} and
eliminate $\SH$. It is easy to check that the resulting coefficients
at $\SJ$ and $\SG$ vanish if either $g = 0$ or
\begin{equation*}
 \eta_{1,1} \, g = -   \theta_2 \xi_0^{-1} \xi_1^{-1} \xi_2^{-1}
 (\xi_0 \xi_1 \xi_2^{-1} + \xi_0 + \xi_1) =
 -  \theta_2^{-1} \xi_0 \xi_1 \xi_2 (\xi_0 + \xi_1) .
\end{equation*}
However, the last equality cannot hold because
$\xi_0^{2} \xi_1^{2} \xi_2^{2} = \theta_2^2$
(see \rf{xixi}).

For $n \,{\geq}\, 3$, the coefficient at $\SJ$ in \rf{FGH}
vanishes only if $g \= 0$.
Thus, the coefficient $g$ in \rf{Rans} must be zero
if $n \geq 2$. Therefore, if $\SR'$ coincides with $\SR$
in the two highest orders, then $\SR' = \SR$.
An analogous statement can be established if we consider
the ansatz \rf{Rans} with $\SR$ being replaced by $\SR^{-1}$.
Thus, Proposition~1 is proven.

\vspace*{2mm}
{\bf Acknowledgments.}
The author thanks P.~Kulish for useful remarks.
This work was supported by the RFBR grants
08--01--00638, 09--01--12150, 09--01--93108.

\appendix

\subsection*{Appendix}

{\em Proof} of Lemma~\ref{PM}.\\
The 6--$j$ symbols of the algebra $U_q(sl_2)$ satisfy
the following q--analogue of the Racah identity \cite{KR1,No}:
\beq{Racah}
\begin{aligned}
 \sum_p {} & \Bigl( (-1)^p \, [2p+1]_q \,
 \FRW{r_1}{r_2}{r_3}{r_4}{l}{p}_{\!q} \,
 q^{\rho(p)-\rho(r_1)-\rho(r_4)}
 \FRW{r_1}{r_3}{r_2}{r_4}{l'}{p}_{\!q} \Bigr) \\
{} & = (-1)^{l+l'} \,
 q^{\rho(r_2)-\rho(l)}
 \FRW{r_3}{r_2}{r_1}{r_4}{l}{l'}_{\!q} \,
 q^{\rho(r_3)-\rho(l')} \,.
\end{aligned}
\end{equation}
(Note that the identity remains true if we set
$\rho(t) = - t(t \+ 1)$, since the 6--$j$ symbols
are self--dual with respect to the replacement $q \to q^{-1}$).

Consider the matrix entry $(kk')$ of equality \rf{AD0}.
Using formula \rf{Dkk} and taking into account that
$A^{(s,n)}$ is a symmetric matrix, we obtain:
\beq{AD0b}
\begin{aligned}
 \sum_m  {}& (-1)^m \, A^{(s,n)}_{k m} \,
 q^{\rho(2s-m)-2\rho(s)} \, A^{(s,n)}_{k' m}  \\
{}&  = (-1)^{n+k+k'} \,
 q^{\rho(3s-n) + \rho(s) -\rho(2s-k) -\rho(2s-k')} \,
 A^{(s,n)}_{k k'} \,.
\end{aligned}
\end{equation}
Now, taking into account formula \rf{A1}, it is easy to see
that relation \rf{AD0b} follows from the identity \rf{Racah}
if we set
$r_1=r_2=r_3=s$, $r_4=3s-n$, $l=2s-k$, $l'=2s-k'$, $p=2s-m$.
\\[-0.5mm]

\noindent
{\em Proof} of Lemma~\ref{PPrel1}. \\
We will prove those relations of Lemma~\ref{PPrel1} that
contain $\SR^{+1}$ on the l.h.s. Their counterparts with
$\SR^{-1}$ on the l.h.s. can be proven similarly.

The second relation in \rf{pp3}:
\begin{eqnarray*}
{}&&
 \pi^{(\ovm)} \hat{\pi}^{(\ovm)} \pi^{(\ovm)} =
 \pi^{(\ovm)} A^{(s,n)} \pi^{(\ovm)} A^{(s,n)} \pi^{(\ovm)} =
 (A^{(s,n)}_{\ovm \ovm})^2 \pi^{(\ovm)} .
\end{eqnarray*}
Here and below we denote
$\hat{\pi}^{(\ovm)} \equiv A^{(s,n)} \pi^{(\ovm)} A^{(s,n)}$.

Relations \rf{pp4}:
\begin{eqnarray*}
{}&&\!\!\!
 \pi^{(\ovm)} \hat{D}_0^{(n)} \pi^{(\ovm)} \rfeq{Dhat}
 \pi^{(\ovm)} A^{(s,n)} D_0^{(n)} A^{(s,n)} \pi^{(\ovm)} \rfeq{AD0}
 \theta_n \pi^{(\ovm)} (D_0^{(n)})^{-1} A^{(s,n)}
    (D_0^{(n)})^{-1} \pi^{(\ovm)} \\
{}&&\!\!\!
 \rfeq{Dkk} \theta_n \xi_{\ovm}^{-2} \pi^{(\ovm)} A^{(s,n)} \pi^{(\ovm)}
 = \theta_n \xi_{\ovm}^{-2} A^{(s,n)}_{\ovm \ovm} \pi^{(\ovm)} ,\\
{}&&\!\!\!
 D_0^{(n)} \hat{\pi}^{(\ovm)} D_0^{(n)} =
 D_0^{(n)} A^{(s,n)} \pi^{(\ovm)} A^{(s,n)} D_0^{(n)} \rfeq{Dkk}
 \xi_{\ovm}^{-2} D_0^{(n)} A^{(s,n)} D_0^{(n)}
    \pi^{(\ovm)} D_0^{(n)} A^{(s,n)} D_0^{(n)} \\
{}&&\!\!\!
 \rfeq{AD0} \theta_n^2 \xi_{\ovm}^{-2}
    A^{(s,n)} (D_0^{(n)})^{-1} A^{(s,n)}
    \pi^{(\ovm)} A^{(s,n)} (D_0^{(n)} )^{-1} A^{(s,n)}
 \rfeq{Dhat} \theta_n^2 \xi_{\ovm}^{-2}  (\hat{D}_0^{(n)})^{-1}
    \pi^{(\ovm)} (\hat{D}_0^{(n)})^{-1} .
\end{eqnarray*}

The first relation in \rf{pp5} (the second can be proven similarly):
\begin{eqnarray*}
{}&&\!\!\!
 \pi^{(\ovm)} \hat{\pi}^{(\ovm)} D_0^{(n)} =
 \pi^{(\ovm)} A^{(s,n)} \pi^{(\ovm)} A^{(s,n)} D_0^{(n)} =
 A^{(s,n)}_{\ovm \ovm} \pi^{(\ovm)}
    A^{(s,n)} D_0^{(n)} (A^{(s,n)})^2 \\
{}&&\!\!\!
 \rfeq{AD0} \theta_n A^{(s,n)}_{\ovm \ovm} \pi^{(\ovm)}
    (D_0^{(n)})^{-1} A^{(s,n)} (D_0^{(n)})^{-1} A^{(s,n)}
 \rfeq{Dkk} \theta_n \xi_{\ovm}^{-1} A^{(s,n)}_{\ovm \ovm}
    \pi^{(\ovm)} (\hat{D}_0^{(n)})^{-1} .
\end{eqnarray*}
\\[-0.5mm]

\noindent
{\em Proof} of Lemma~\ref{LD}.\\
For $n \=1$, the matrices $\SG$, $\SH$, $\SJ$ are
of the size $2{\times}2$ and relations \rf{JGH1}
can be verified straightforwardly using the explicit
form of the matrix  $A^{(s,1)}$ (see eq. (73) in \cite{B1}).

In order to examine the case $n \,{\geq}\,2$, let us
write down explicitly the matrix entries of $\SG$, $\SH$, and $\SJ$:
\begin{align}
\label{Gkk}
 \SG_{kk'} &= \delta_{kn} \, \delta_{k' n} \,  -
    A^{(s,n)}_{nk}  A^{(s,n)}_{nk'} ,\\
\label{Hkk}
 \SH_{kk'} &= \theta_n^{-1} \xi_n
 ( \delta_{kn} \, \xi_{k'} A^{(s,n)}_{nk'} +
 \delta_{k'n} \, \xi_{k} A^{(s,n)}_{nk} ) -
 (\xi_{k}^{-1} \+ \xi_{k'}^{-1})\, A^{(s,n)}_{nk} A^{(s,n)}_{nk'} , \\
\label{Jkk}
 \SJ_{kk'} &= (\xi_{k} \xi_{k'} -
    \theta_n^2 \xi_n^{-2} \xi_{k}^{-1} \xi_{k'}^{-1} )\,
    A^{(s,n)}_{nk} A^{(s,n)}_{nk'} .
\end{align}
Recall that $k,k' \=0,1,{\ldots},n$.

Considering \rf{Jkk} for $k \= 0$ and $k' \= 0,1$, it is
easy to infer that $\SJ \,{\neq}\, 0$ (since $\xi_0^2 \,{\neq}\, \xi_1^2$).

Assume that the following relation holds
\beq{abg}
\alpha \SG  + \beta \SJ - \gamma \SH  =0\,,
\end{equation}
where $\alpha \beta \gamma \,{\neq}\, 0$.
Using formulae \rf{Gkk}--\rf{Jkk}, write down the matrix
entries of \rf{abg} for $(k,k') \=(0,0)$,
$(0,1)$, $(1,1)$ dividing them by $A^{(s,n)}_{nk} A^{(s,n)}_{nk'}$
(note that $A^{(s,n)}_{nk} \,{\neq}\,0$ for all $k$,
see eq. (97) in \cite{B1}):
\begin{align}
\nonumber
 -\alpha + (\xi_0^2 - \theta_n^2 \xi_n^{-2} \xi_0^{-2}) \beta
 + 2 \xi_0^{-1} \gamma = 0 \,, \\
\label{abg'}
 -\alpha + (\xi_0 \xi_1 -
    \theta_n^2 \xi_n^{-2} \xi_0^{-1} \xi_1^{-1}) \beta
 + (\xi_0^{-1} + \xi_1^{-1}) \gamma = 0 \,, \\
\nonumber
 -\alpha + (\xi_1^2 - \theta_n^2 \xi_n^{-2} \xi_1^{-2}) \beta
 + 2 \xi_1^{-1} \gamma = 0 \,.
\end{align}
The determinant of this system of equations is
$d \= (\xi_0^{-1}{-}\xi_1^{-1})^3
 (\theta_n^2 \xi_n^{-2}{-} \xi_0^{2}\xi_1^{2})$.
Since $\xi_0 \,{\neq}\, \xi_1$ then the equality
$d \=0$ can be satisfied only if
\beq{xixi}
  \theta_n^2 =  \xi_0^{2}\xi_1^{2} \xi_n^{2} \,,
\end{equation}
which is equivalent to the following condition:
$\rho(3s{-}n) \+ 3\rho(s) \,{-}\, \rho(2s) -
    \rho(2s{-}1) \,{-}\, \rho(2s{-}n) \= 2s(2 \,{-}\,n) \= 0$.
Thus, relation \rf{abg} cannot hold for $n \,{\geq}\,3$.

For $n \=2$, a solution of the system \rf{abg'} is given by:
$\alpha \= \beta^{-1} \= \xi_0 \+ \xi_1$,
$\gamma \= \xi_0 \xi_1$.
A direct check, using the explicit form of the matrix
$A^{(s,2)}$ (see eq. (74) in \cite{B1}), shows that relation
\rf{abg} with such coefficients holds indeed.
Since system \rf{abg'} has no solution for $\gamma \= 0$,
we conclude that $\SG$ and $\SJ$ are linearly independent.

%
\small  \setlength{\itemsep}{-3pt}
\newcommand{\myitem}[6]{#1:\, {\em #2}.\,--- {#3} {\bf #4} (#5), #6.}

\end{document}